\documentclass[12pt]{article}
\author{Laurent Bartholdi\footnote{Department of Mathematics,
University of California, Berkeley, CA 94720. {\sf
laurent@math.berkeley.edu.}} \ and B\'alint
Vir\'ag\footnote{Department of Mathematics, MIT, Cambridge, MA
02139. {\sf balint@math.mit.edu}. Research partially supported
by NSF grant \#DMS-0206781.}}
\title{Amenability via random walks}
\date{May 19, 2003}
\def\myabstract{We answer an open question of Grigorchuk and
{\.Z}uk about amenability using random walks. Our results
separate the class of amenable groups from the closure of
subexponentially growing groups under the operations of group
extension and direct limits; these classes are separated even
within the realm of finitely presented groups.}

\newtheorem{theorem}{Theorem}

\newtheorem{proposition}[theorem]{Proposition}
\newtheorem{lemma}[theorem]{Lemma}
\newtheorem{corollary}[theorem]{Corollary}
\newenvironment{theorem*}[1]{\par \trivlist
 \itemindent 0pt \item[\hskip\labelsep\bf Theorem #1]
 \it\ignorespaces}{\endtrivlist}

\newcommand{\qed}{\hfill\mbox{$\framebox(5,5)[]{}$}}
\newenvironment{proof}{\par \trivlist
 \itemindent\parindent \item[\hskip\labelsep\sc Proof.]
 \ignorespaces}{\qed\endtrivlist}

\DeclareSymbolFont{AMSb}{U}{msb}{m}{n}
\DeclareSymbolFontAlphabet{\mathbb}{AMSb}

\newcommand\mnote[1]{}
\newcommand\lb[1]{\label{#1}\mnote{#1}}
\newcommand\bel[1]{{\mnote{#1}}\begin{equation}\label{#1}}
\newcommand\ee{\end{equation}}

\newcommand{\square}{{\framebox(2,2)[]{}}}
\newcommand\re[1]{(\ref{#1})}

\newcommand{\eps}{\varepsilon}

\newcommand{\ev}{{\mathbf E}}
\newcommand{\pr}{{\mathbf P}}

\newcommand{\as}{\mbox{\hspace{.3cm} a.s.}}

\newcommand{\AG}{\mathsf{AG}}
\newcommand{\EG}{\mathsf{EG}}
\newcommand{\SG}{\mathsf{SG}}

\usepackage{graphicx}
\usepackage{natbib}
\begin{document}
\maketitle \begin{abstract}{\noindent \myabstract}\end{abstract}
\section{Introduction}

The concept of amenability, introduced by \cite{vonneumann},
has been central to many areas of mathematics. \cite{kesten}
showed that a countable group is amenable if and only if the
spectral radius equals 1; in particular, if the random walk
escapes at a sublinear rate. Although this connection has been
deeply exploited to study the properties of random walks, it
appears that it has not yet been used to prove the amenability
of groups.

A group is amenable if it admits a finitely additive invariant
probability measure. The simplest examples of amenable groups
($\AG$)  are
{\list{}{\itemsep 0in \topsep 1ex}
 \item[(i)] finite and Abelian groups and, more generally,
 \item[(ii)]  groups of subexponential growth.
 \endlist}
Amenability is preserved by taking subgroups, quotients,
extensions, and direct limits. The classes of elementary amenable
($\EG$), and subexponentially amenable ($\SG$, see \cite{g98}, and
\cite{cgh}, \S14) groups are the closure of (i), (ii) under these
operations, respectively. We have
$$
\EG\subseteq \SG\subseteq \AG,
$$
and the question arises whether these inclusions are strict:
 \cite{day} asked this about $\EG \subseteq
\AG$ (see also \cite{g98}). \cite{chou} showed that there are
no elementary amenable groups of intermediate growth. Thus
Grigorchuk's group separates the class $\EG$ and $\SG$,
answering Day's question.

In this paper, we show by example that the inclusion
$\SG\subset \AG$ is also strict.

The group $G$ we are considering is the iterated monodromy
group of the polynomial $z^2-1$. It was first studied by
\cite{gz02a}, who showed  that $G$ does not belong to the class
$\SG$. The main goal of this note is to show, using rate of
escape for random walks, that $G$ is amenable. This answers a
question of the above authors.

\begin{figure}
\label{julia} \centering
\includegraphics[height=2.2in]{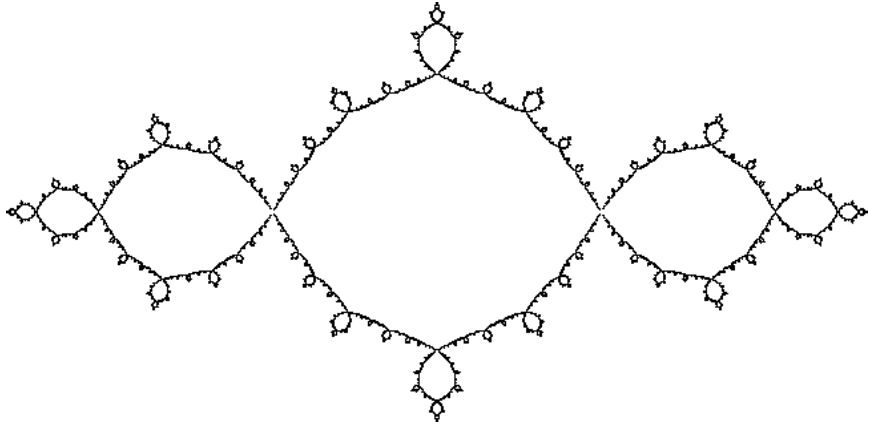}
\caption{The basilica, or the action of $G$ on $\partial T$}
\end{figure}

Let $T$ be the rooted binary tree with vertex set $V$
consisting of all finite binary sequences, and edge set $
E=\{(v,vi)\,:\, v\in V,\, i\in \{0,1\}\}. $ Let
$\eps\in\mbox{Aut}(T)$ send $iv$ to $((i+1)$ mod $2)v$. For
$g,h\in$Aut$(T)$ (with the notation $g:v\mapsto v^g$) let
$(g,h)$ denote the element of Aut$(T)$ sending $0v\mapsto 0
v^g$ and $1v \mapsto 1 v^h$. The group $G$ is generated by the
following two recursively defined elements of Aut$(T)$:
$$
a=(1,b),\quad b=(1,a)\eps.
$$
Then $G$ is the iterated monodromy group of the polynomial
$z^2-1$; the scaling limit of the Schreier graphs of its action
on level $n$ of $T$ is the {\it basilica}, i.e. the Julia set
of this polynomial (see the survey \citet*{bgn} and Figure
\ref{julia}).
Let $Z_n$ be a symmetric random walk on $G$ with step
distribution supported on $a,\,b$ and their inverses, and let
$|\cdot|$ denote shortest-word distance.

\begin{theorem}
We have
$\lim |Z_n|/n = 0$ with probability 1.
\end{theorem}

This implies that the heat kernel decays sub-exponentially, the
group $G$ is amenable (see \cite{kw}), and $\SG\not=\AG$.

In the rest of the paper, we extend this result in two
directions. In Section \ref{s.fg}, Theorem \ref{finitely} we
give a finitely presented example $\tilde G$  separating $\AG$
and $\SG$, showing that these classes are distinct even in this
realm. \cite{g98} showed  that $\EG\not= \AG$ (more precisely,
$\EG \not= \SG$) for finitely presented groups.

 In Section \ref{s.ss} Corollary \ref{5/6}, we give
a quantitative upper bound of order $n^{5/6}$ on the rate of
escape. For the heat kernel, we have the following quantitative
lower bound.

\begin{theorem}
There exists $c>0$ so that for all $n$ we have
$\pr(Z_{2n}=1)\ge e^{-cn^{2/3}}.$
\end{theorem}

Motivated by their question about amenability, \cite{gz02b}
study spectral properties of $G$. Amenability of $G$  has been
claimed in the preprint \cite{barth}, whose proof appears to
contain serious gaps and is considered altogether incomplete.
The present paper uses the same starting point as \cite{barth},
but follows a different path; we get specific heat kernel
bounds for a less general family of groups.

\section{A fractal distance}

For $g\in G$ and $v\in T$ let $g[v]\in$Aut$(T)$ denote the action
of $g$ on the descendant subtree of $v$, and let $g(v)\in C_2$
denote the action on the two children of $v$. Let $S$ be a finite
binary subtree of $T$ containing the root (i.e. each vertex in $S$
has zero or two descendants). Let $\partial S$ denote its set of
leaves, and $E=\# \partial S-1$ the number of edges. Let
 \begin{eqnarray*}
 \nu_S(g)= E + \sum_{v\in \partial S} |g[v]|,\\
 \nu(g) = \min_S {\nu_S(g)}.
 \end{eqnarray*}
The quantity $\nu$ has the alternative recursive definition; for
$g=(g_1,g_2)\eps_*$, let
 $$
 \nu(g)=\min(|g|,1+\nu(g_1)+\nu(g_2)).
 $$

\begin{lemma} The function $\nu$ is a norm on $G$. Moreover, $\nu$-balls have
exponential growth.
\end{lemma}
\begin{proof}
First note that since multiplying $g$ by a increases
$|g_1|+|g_2|$ by at most $1$, we get $ |g|\ge |g_1|+|g_2| $.
This implies that if $\nu(g)=|g|$ then $\nu(g)\ge
\nu(g_1)+\nu(g_2)$. So in general, we have
 \bel{fractaldist}
 \nu(g_1)+\nu(g_2) \le \nu(g) \le
\nu(g_1)+\nu(g_2)+1.
 \ee
We now check that $\nu$ satisfies the triangle inequality; this is
clear if $\nu(g)=|g|,\,\nu(h)=|h|$, otherwise we may assume that
$\nu(g)=\nu(g_1)+\nu(g_2)-1$. Then we get
$$
\nu(gh)\le 1+\nu((gh)_1)+\nu((gh)_2) \le
1+\nu(g_1)+\nu(g_2)+\nu(h_1)+\nu(h_2) \le \nu(g)+\nu(h),
$$
where the first inequality holds by induction (some care is needed
to show that the induction can be started).

We now claim that the balls $B_n=\{g\,:\, \nu(g)\le n\} $ grow
at most exponentially, more precisely, we have
\bel{balls} \# B_n \le
40^n   \quad \mbox{for all } n.
\ee
Indeed, there are at most $4^n$ such subtrees $S$ with at most $n$
edges. Given the subtree $S$, the element $g\in B_n$ is defined by
its action $g(v)\in C_2$ at the vertices of $S$ that are not
leaves (at most $2^n$ possibilities), as well as the words $g[v]$
at the vertices $v$ that are leaves (these can be described with
$n$ symbols from the alphabet $a,a^{-1},b,b^{-1}$ and comma). Thus
we have $\#B_n \le (4\cdot 2 \cdot 5)^n$.

For the other direction, note that $\nu$-balls contain the
word-distance balls of the same radius and  $G$ has exponential
growth (see \cite{gz02a}).
\end{proof}

\section{Self-similarity of random walks on $G$}\label{s.ss}

Fix $r>0$, and consider the random walk $Z_n$ on the free group
$F_2$ where each step is chosen from $(a,a^{-1},b,b^{-1})$
according to weights $(1,1,r,r)$, respectively. This walk projects
to $F_2 \wr C_2$ via the substitution $a\mapsto (1,b)$ and
$b\mapsto(1,a)\eps$. Let $(Y_n,X_n)\eps_n$ be the projection of
$Z_n$. Define the stopping times
\begin{eqnarray*}
\sigma(0)&=&0, \\
\sigma(m+1)&=&\min\{n>\sigma(m)\,:\, \eps_n=1,\,X_n\not=X_{\sigma(m)}\}, \quad m\ge 0,\\
\tau(0)&=&\min\{n>0\,:\,\eps_n=\eps\},\\
\tau(m+1)&=&\min\{n>\tau(m)\,:\, \eps_n=\eps,\,Y_n\not=Y_{\tau(m)}\}, \quad m\ge 0.
\end{eqnarray*}

\begin{lemma}\lb{fractalwalk}
$X_{\sigma(m)}$, $Y_{\tau(m)}$ are simple random walks on $F_2$
with step distribution given by the weights $(r/2,r/2,1,1)$ on
$(a,a^{-1},b,b^{-1})$, respectively.
\end{lemma}
\begin{figure}[t]
  \begin{center}
 \includegraphics[height=2.5in]{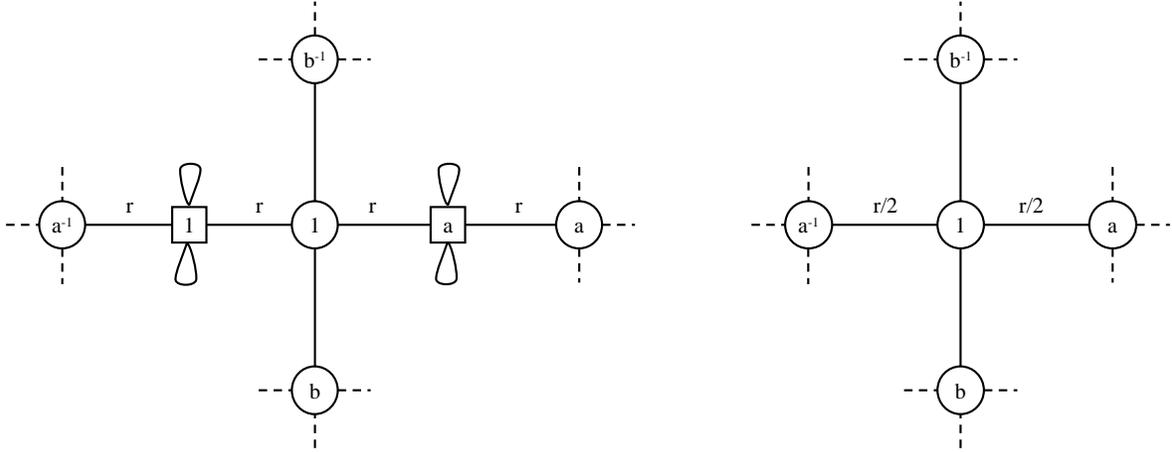}
    \caption{The random walks $(X_n,\eps_n)$ and $X_{\sigma(m)}$}
    \label{fig:pif1}
  \end{center}
\end{figure}
\begin{proof}
The process $(X_n,\eps_n)$ (i.e. we ignore the rest of the
information in $Z_n$) is a random walk on a weighted graph, as
shown on the left of Figure \ref{fig:pif1}; there the value of
$\eps_n$ is represented by a circle $(1)$ or square $(\eps)$.
 When
we only look at this walk at the times $\sigma(n)$, the resulting
process is still a Markov chain, where the transition
probabilities are given by the hitting distribution on the $4$
circles that are neighbors or separated by a single square.

The hitting distribution is given by effective conductances,
and using the series law we get the picture on the right hand
side of Figure \ref{fig:pif1}. This process is a symmetric
random walk with weights as claimed.

The proof for $Y$ is identical. Because $\tau(0)$ has a different
definition, the process $Y_{\tau(m)}$ does not start at $0$,
rather at $a$ or $a^{-1}$. Note that the processes
$X_{\sigma(m)}$ and $Y_{\tau(m)}$ are not
independent.
\end{proof}

\begin{lemma}\label{sigma} With probability 1,  we have
$ \lim m/\sigma(m)=(2+r)/(4+4r)=:f(r)$, and the same holds for $\tau$.
\end{lemma}

\begin{proof}
The increments of the process $\sigma(m)$ are the time the random
walk in Figure \ref{fig:pif1} spends between hitting two different circles.
These increments are independent and identically distributed. Let
$t_\circ,\, t_\square$ denote the expected times starting from a
circle or a neighboring square to hit a different circle.
Conditioning on the first step of the walk gives the equations
\begin{eqnarray*}
t_\circ&=& 1+r/(r+1)t_\square \\
t_\square &=& 1+r/(2(r+1)) t_\circ + 1/(r+1) t_\square
\end{eqnarray*}
And the solution is $t_\circ=4(1+r)/(2+r)$. The claim now follows
from the strong law of large numbers.
\end{proof}

If $\overline Z_n$ denotes the image of $Z_n$ in $G$,
then by construction we have $\overline Z_n=(\overline Y_n,
\overline X_n)\eps_n$. In the rest of this section we will
simply (ab)use the notation $Z_n,\, X_n,\, Y_n$ for the images in $G$ of
the corresponding random words.

\begin{proposition} \lb{nuspeed} We have $\lim \nu(Z_n)/n = 0 \as$.
\end{proposition}

\begin{proof}
By Kingman's subadditive ergodic theorem, the random limit
$$
s(r)=\lim \nu(Z_n)/n
$$
exists and equals a constant with probability 1. By
(\ref{fractaldist}) we also have
$$
\nu(Z_n) \le \nu(X_n) + \nu(Y_n) +1
$$
and therefore
\begin{eqnarray*}
s(r) &\le&  \limsup \nu(X_n)/n + \limsup
\nu(Y_n)/n  \\
&=& \limsup \nu(X_{\sigma(m)})/\sigma(m) + \limsup
\nu(Y_{\tau(m)})/\tau(m)
\\
&=& \left(\lim \nu(X_{\sigma(m)})/m\right) \left(\lim
m/\sigma(m)\right) + \left(\lim \nu(Y_{\tau(m)})/m\right)
\left(\lim m/\tau(m)\right) \\
&=& 2s(2/r)f(r).
\end{eqnarray*}
In the last equality we used Lemmas \ref{fractalwalk} and
\ref{sigma}. Iterating this inequality and we get
$$
s(r) \le 4s(r)f(r)f(2/r)=4s(r)/8
$$
and since $s$ is a finite constant we get $s=0$ with probability 1.
\end{proof}

Proposition \ref{nuspeed} implies that the walk has asymptotic
entropy 0 (see the proof of Theorem 5.3 in KW (2002)). This, by
Theorem 5.13 ibid. implies the following.

\begin{corollary}
The random walk on $G$ satisfies $\lim |Z_n|/n=0$ a.s.
\end{corollary}

A more technical version of Proposition \ref{nuspeed} gives a
better bound on the $\nu$-rate of escape.

\begin{proposition}\label{tech} There exists $c>0$ depending on $r$ so that for all
$n\ge 1$
$$u_r(n):=\ev \max_{i\le n} \nu(Z_i) \le cn^{2/3}.$$
\end{proposition}

\begin{proof}
Let $L(n)$ be the largest so that $\sigma(L)\le n$, and let $M(n)$
be the largest so that $\tau(M)\le n$. Following the argument of
Proposition \ref{nuspeed}, we get
\begin{eqnarray*}
u_r(n) &\le&
\vphantom{\frac{1}{2}}\ev \max_{i\le n} \nu(X_i) + \ev \max_{i\le n}\nu(Y_i) +1 \\
&=& \vphantom{\frac{1}{2}}\ev \max_{i\le L(n)}\nu(X_{\sigma(i)})+
\ev \max_{i\le M(n)} \nu (Y_{\tau(i)}) + 1\\
&\le &  \vphantom{\frac{1}{2}}2 + 2 u_{2/r}(f(r)n+k_n) +
2n\pr(L(n)>f(r)n+k_n)
\end{eqnarray*}
where we can choose the constants $k_n=n^{\alpha}$ for some
$\alpha>1/2$ close to $1/2$. By the large deviation principle and
the fact that the distribution of $\sigma(1)$ has an exponential
tail, the last term can be bounded above by
$c_1ne^{-c_2n^{2\alpha-1}}<c_3$. Thus for $n\ge 1$ we get
$$
u_r(n) \le u_{2/r}(f(r)n+k_n) + c_4 \le u_{2/r}(f(r)n) + (c_4+1) n^\alpha.
$$
Applying this to $u_{2/r}$ as well and using the fact that
$f(r)f(2/r)=1/8$, we easily get
$$
u_r(8n) \le 4u_r(n)+ c_6 n^{\alpha}
$$
with $u_r(1)\le 1$. Iteration at values of $n$ that are powers of
$8$ gives that for such values
$$
u_r(n) \le c_7 n^{2/3}.
$$
Since $u_r(n)$ is monotone in $n$, the claim follows.
\end{proof}

\begin{corollary}
There exists $c>0$ depending on $r$ so that for all $n\ge 0$ we have
$$\pr(Z_{2n}=1)\ge e^{-cn^{2/3}}.$$
\end{corollary}

\begin{proof}
By Markov's inequality, we have
$$
\pr(\nu(Z_n) \le 2cn^{2/3}) \ge 1/2
$$
and therefore there exists $g\in B_{2cn^{2/3}}$ so that
$$\pr(Z_n=g)\ge 1/(2\#B_{2cn^{2/3}})$$
and since balls grow at most exponentially (\ref{balls}), we get
$$
\pr(Z_{2n}=1) \ge \pr(Z_n=g)\pr(Z_n^{-1}Z_{2n}=g^{-1})
=\pr(Z_n=g)^2 \ge c_1e^{-c_2n^{2/3}}
$$
and the claim follows.
\end{proof}

Let $M_n=\max(|X_1|,\ldots ,|X_n|)$. We have the following
bound on the rate of escape.

\begin{corollary}\label{5/6}
There exists $c>0$ so that for all $a,n\ge 1$ we have $\pr(M_n>
an^{5/6}) < c/a$.
\end{corollary}

\begin{proof}
Let $K_n=\max(\nu(X_1),\ldots ,\nu(X_n))$. We have
 \begin{eqnarray}\nonumber
\pr(M_n > an^{5/6})&=&\pr(M_n > an^{5/6},K_n> ac_1n^{2/3})\\&+&
\pr(M_n> an^{5/6},K_n\le ac_1n^{2/3}).\label{twoterms}
 \end{eqnarray}
The first term is at most $\pr(M_n> ac_1n^{2/3})<c_2/(c_1a)$ by
Proposition \ref{tech} and Markov's inequality. The second term
is bounded above by the sum of $\pr(X_m=g)$ over all $m\le n$
and all $g$ with $|g|>an^{5/6}$ and  $\nu(g)\le ac_1n^{2/3}$.
By the Varopoulos-Carne bounds (see KW, 2002) the first
constraint on $g$ implies
$$
\pr(X_m=g) \le e^{-(an^{5/6})^2/(2n)},
$$
and since $\nu$-balls grow exponentially (\ref{balls}) the
second term of \re{twoterms} is bounded above by
$$
n e^{
-(an^{5/6})^2/(2n)+c_3c_1an^{2/3}}=ne^{(c_3c_1a-a^2/2)n^{2/3}},
$$
which is at most $c_4/a$ for an appropriate choice of $c_1$.
\end{proof}

\section{A finitely presented example and
generalizations}\label{s.fg}

Our first goal is to show that an \textsf{HNN}-extension of $G$
gives a finitely presented example separating $\AG$ and $\SG$.
The following lemma is needed.

\begin{lemma}\label{thm:presG}
  $G$ has the following presentation:
  \[G = \langle a,b\,|\,\sigma^n[a,a^b]\quad\forall n\in\mathbb{N}\rangle,\]
  where $\sigma$ is the substitution $b\mapsto a,a\mapsto b^2$.
\end{lemma}

\begin{proof}
  By Proposition 9 of \cite{gz02a} we have
  \[G = \langle a,b\,|\,\sigma^n[a,a^{b^{2m+1}}]
     \quad\forall n,m\in\mathbb{N}\rangle.\]
  For odd $i$, we have $a^{b^i}\equiv[a^{-1},b^{-2}]^ba^{b^{i-2}}$
  using the relation $[b^{2a},b^2]=\sigma([a^b,a])$; therefore
  $[a^{b^i},a]$ follows from $[a^{b^{i-2}},a]$ and
  $[[a^{-1},b^{-2}]^b,a]$, which itself is a consequence of $[a^b,a]$.
  So the relations $\sigma^n([a,a^{b^{2m+1}}])$ may be
  eliminated for all $m>1$ as long as $\sigma^n([a,a^b])$ and
  $\sigma^{n+1}([a,a^b])$ are kept.
\end{proof}

\begin{theorem}\label{finitely}
  $G$ embeds in the finitely presented group
  \[\tilde G = \langle a,t\,|\,a^{t^2}=a^2,\,[[[a,t^{-1}],a],a]=1\rangle.\]
  Furthermore, $\tilde G$ is also amenable, and does not belong to the
  class $\SG$.
\end{theorem}
This implies that the classes \textsf{SG} and \textsf{AG} are
distinct, even in the realm of finitely presented groups.

\begin{proof}
  Let $\tilde G$ be the \textsf{HNN} extension of $G$ along the
  endomorphism $\sigma$ identifying $G$ and $\sigma(G)$: it is given
  by the presentation
  \[\tilde G=\langle
  a,b,t|\,a^t=\sigma(a),\,b^t=\sigma(b),\,\mbox{relations in }G\rangle.\]
  A simpler presentation follows by eliminating the generator $b$.

  Consider the kernel $H$ of the map $a\mapsto 1,t\mapsto t$ from
  $\tilde G$ to $\langle t\rangle$. Since the \textsf{HNN} extension
  is ``ascending'', we have $H=\bigcup_{n\in\mathbb{Z}} G^{t^n}$, an ascending
  union. Therefore $H$ is amenable, and since $\tilde G$ is an
  extension of $H$ by $\mathbb Z$, it is also amenable.

  Finally, if $\tilde G$ were in \textsf{SG}, then $G$ would also be
  in \textsf{SG}, since it is the subgroup of $\tilde G$
  generated by $a$ and $a^{t^{-1}}$. However, Proposition 13 of
  \cite{gz02a} shows that $G$ is not in $\SG$.
\end{proof}

{\bf Generalizations.} In what setting does the proof for
amenability work? Let $G$ be a group acting spherically
transitively on a $b$-ary rooted tree ($b\ge 2$), and suppose
that it is defined recursively by the set $S$ of generators
$g_i=(g_{i,1},\ldots ,g_{i,b})\sigma_i$, where each $g_{i,v}$
is one of the $g_j$. Consider the Schreier graph of the action
of $G$ on $T_1$, that is level $1$ of the tree; we label level
$1$ of the tree by the integers $1,\ldots b$. Furthermore, we
label each directed edge $(v,g_i)$ by $g_{i,v}$.

Fix a vertex at level $1$, without loss of generality the vertex
$1$. Consider the set of cycles in the Schreier graph that go from
1 to 1 and may traverse edges either forwards or backwards; such a
cycle is called ``irreducible'' if it only visits 1 at its
endpoints. The label of a cycle is the product of the labels along

If $\#S\ge 2$ then there are infinitely many irreducible cycles. A
necessary condition (1) for our proof to work is that the set of
labels of irreducible cycles is finite and agrees with the set $S$
of generators together with the identity.

Given a probability distribution $\mu$ on the set of
generators, we get a distribution on the set of irreducible
loops by considering the path of a random walk on $G$ up to the
first positive time $\tau$ that it fixes $1$ and has a cycle
whose label is not $1$. Call the distribution of this label
$\mu'$. The transformation $\mu\mapsto \mu'$ is a continuous
map from a convex set to itself, so it has a fixed point.

A further necessary condition is that at least one fixed point
is in the interior of the convex set, i.e. assigns positive
weight to each generator. For this, it is sufficient that
condition (1) does not hold for any proper subset of $S$.

Now let $\mu_0$ be a such a fixed point, and let $\alpha=\log
b/\log \ev \tau$ for the corresponding random time $\tau$. If
$\alpha>1/2$, then the argument above gives a heat kernel lower
bound of $e^{-cn^\alpha}$.  The argument above cannot give an
exponent below $1/2$ as the rate of escape cannot be slower
than $n^{1/2}$. In the proof, the large deviation bounds for
$\sigma$ break down at $\alpha=1/2$.
\bigskip

\noindent {\bf Example.} Consider the group acting on the
binary tree generated by $a_i=(1,a_{i+1})$ for $i<k$, and
$a_k=(1,a_1)\eps$. The distribution $\mu=(m_1,\ldots, m_k)$ on
the generators (and symmetrically on their inverses) is then
sent to $T \mu'=(m_k/2,m_1,\ldots,m_{k-1})/(1-m_k/2)$. A fixed
point is given by $(1,2^{1/k},\ldots,2^{(k-1)/k})$ normalized
to be a probability distribution. A simple computation gives
$\ev \tau=2^{1+1/k}$, and we get the heat kernel lower bound
$e^{-c n^\alpha}$ with $\alpha=k/(k+1)$.

In this example, it is not important to consider a fixed point.
Since $T^k=1$, one may iterate the decomposition process $k$
times starting from an arbitrary $\mu$. Then one is lead to
consider $2^k$ processes having the same distribution as the
original walk, each with time running slower by a constant
factor. After massive cancellations, one finds that the
constant $\prod \ev \tau_i$ does not depend on $\mu$, and
equals $2^{k+1}$. This gives the same heat kernel bound as
above.
\bigskip

\noindent {\bf Acknowledgments.} B.V. thanks Mikl\'os Ab\'ert for
inspiring discussions and for bringing this problem to his
attention.

\bibliography{pi}

\begin{thebibliography}{11}
\expandafter\ifx\csname natexlab\endcsname\relax\def\natexlab#1{#1}\fi
\expandafter\ifx\csname url\endcsname\relax
  \def\url#1{{\tt #1}}\fi

\bibitem[Bartholdi(2002)]{barth}
L.~Bartholdi (2002).
\newblock {Amenability of groups acting on trees}.
\newblock { arXiv:math.GR/0204076}.
\newblock Version 4, July 13, 2002.

\bibitem[Bartholdi et~al.(2003)Bartholdi, Grigorchuk, and Nekrashevych]{bgn}
L.~Bartholdi, R.~I. Grigorchuk, and V.~V. Nekrashevych.
\newblock From fractal groups to fractal sets.
\newblock In P.~Grabner and W.~W. T.~U. Graz), editors, {\em Fractals in Graz},
  Trends in Mathematics, pages 25--118. Birkha\"user Verlag, Basel, 2003.

\bibitem[Ceccherini et~al.(1999)Ceccherini, Grigorchuk, and de~la Harpe]{cgh}
T.~G. Ceccherini, R.~I. Grigorchuk, and P.~de~la Harpe (1999).
\newblock Amenability and paradoxical decompositions for pseudogroups and
  discrete metric spaces.
\newblock {\em Tr. Mat. Inst. Steklova}, 224\penalty0 (Algebra. Topol. Differ.
  Uravn. i ikh Prilozh.):\penalty0 68--111.

\bibitem[Chou(1980)]{chou}
C.~Chou (1980).
\newblock Elementary amenable groups.
\newblock {\em Illinois J. Math.}, 24\penalty0 (3):\penalty0 396--407.

\bibitem[Day(1957)]{day}
M.~M. Day (1957).
\newblock Amenable semigroups.
\newblock {\em Illinois J. Math.}, 1:\penalty0 509--544.

\bibitem[Grigorchuk(1998)]{g98}
R.~I. Grigorchuk (1998).
\newblock An example of a finitely presented amenable group that does not
  belong to the class {EG}.
\newblock {\em Mat. Sb.}, 189\penalty0 (1):\penalty0 79--100.

\bibitem[Grigorchuk and {\.Z}uk(2002{\natexlab{a}})]{gz02a}
R.~I. Grigorchuk and A.~{\.Z}uk (2002{\natexlab{a}}).
\newblock On a torsion-free weakly branch group defined by a three state
  automaton.
\newblock {\em Internat. J. Algebra Comput.}, 12\penalty0 (1-2):\penalty0
  223--246.
\newblock International Conference on Geometric and Combinatorial Methods in
  Group Theory and Semigroup Theory (Lincoln, NE, 2000).

\bibitem[Grigorchuk and {\.Z}uk(2002{\natexlab{b}})]{gz02b}
R.~I. Grigorchuk and A.~{\.Z}uk.
\newblock Spectral properties of a torsion-free weakly branch group defined by
  a three state automaton.
\newblock In {\em Computational and statistical group theory (Las Vegas,
  NV/Hoboken, NJ, 2001)}, volume 298 of {\em Contemp. Math.}, pages 57--82.
  Amer. Math. Soc., Providence, RI, 2002{\natexlab{b}}.

\bibitem[Kaimanovich and Woess(2002)]{kw}
V.~A. Kaimanovich and W.~Woess (2002).
\newblock Boundary and entropy of space homogeneous {M}arkov chains.
\newblock {\em Ann. Probab.}, 30\penalty0 (1):\penalty0 323--363.

\bibitem[Kesten(1959)]{kesten}
H.~Kesten (1959).
\newblock Full {B}anach mean values on countable groups.
\newblock {\em Math. Scand.}, 7:\penalty0 146--156.

\bibitem[von Neumann(1929)]{vonneumann}
J.~von Neumann (1929).
\newblock Zur allgemeinen {Theorie} des {Masses}.
\newblock {\em Fund. Math.}, 13:\penalty0 73--116 and 333.
\newblock = \emph{Collected works}, vol.\ I, pages 599--643.

\end{thebibliography}

\end{document}